\magnification=1200
\overfullrule=0pt
\noindent
{\bf Infinitely many local minima of sequentially weakly lower semicontinuous
functionals and applications}\par
\bigskip
\bigskip
{\centerline {BIAGIO RICCERI}}\par
\bigskip
\bigskip
\bigskip
\bigskip
I am here concerned with an overview of some applications
of Theorem 1 below that have been obtained in the last two years.\par
\bigskip
THEOREM 1 ([14]). - {\it Let $X$ be a non-empty sequentially weakly closed
set in a reflexive
 real Banach space, and let $\Phi, \Psi
:X\to ]-\infty,+\infty]$ be two sequentially weakly lower semicontinuous
functionals.
 Assume also that $\Psi$ is
(strongly) continuous. Denote by $I$ the set of all real numbers $\rho>
\inf_{X}\Psi$ such that $\Psi^{-1}(]-\infty,\rho[)$ is bounded and
intersects the domain of $\Phi$. Assume that $I\neq \emptyset$.
 For each $\rho\in I$, put 
$$\varphi(\rho)=\inf_{x\in \Psi^{-1}(]-\infty,\rho[)}
{{\Phi(x)-\inf_{\overline {(\Psi^{-1}(]-\infty,\rho[))}_{w}}\Phi}
\over {\rho-\Psi(x)}}\ ,$$
where $\overline {(\Psi^{-1}(]-\infty,\rho[))}_{w}$ is the closure
of $\Psi^{-1}(]-\infty,\rho[)$ in the relative weak topology of $X$.
 Furthermore, set
$$\gamma=\liminf_{\rho\to (\sup I)^{-}}\varphi(\rho)$$
and
$$\delta=\liminf_{\rho\to 
\left ( {\inf\limits _{X}\Psi}\right ) ^{+}}\varphi(\rho)\ .$$
Then, the following conclusions hold:\par
\smallskip
\noindent
$(a)$\hskip 10pt For each $\rho\in I$ and each $\mu>
\varphi(\rho)$, the restriction of the
functional $\Phi+\mu\Psi$ to $\Psi^{-1}(]-\infty,\rho[)$ has a global
minimum.\par
\smallskip
\noindent
$(b)$\hskip 10pt If $\gamma<+\infty$, then, for each $\mu>
\gamma$, the following alternative holds: either the restriction of
$\Phi+\mu\Psi$ to $\Psi^{-1}(]-\infty,\sup I[)$ has a global minimum,
 or there exists a sequence $\{x_{n}\}$
of local minima of $\Phi+\mu\Psi$ such that $\Psi(x_{n})<\sup I$ for
all $n\in {\bf N}$ and
$\lim_{n\to \infty}\Psi(x_{n})=\sup I$.\par
\smallskip
\noindent
$(c)$\hskip 10pt If $\delta<+\infty$, then, for each $\mu>
\delta$, 
 there exists a
sequence $\{x_{n}\}$ of local minima of
$\Phi+\mu\Psi$, with $\lim_{n\to \infty}\Psi(x_{n})=\inf_{X}\Psi$,
 which weakly converges to a global minimum of
$\Psi$.}\par
\bigskip
A first consequence of conclusion $(a)$ of Theorem 1 is as follows:\par
\bigskip
THEOREM 2 ([14]). - {\it Let $E$ be a reflexive real Banach space, $X$ a closed,
convex, unbounded subset of $E$, and $\Phi, \Psi :X\to {\bf R}$ two convex
functionals, with $\Phi$ lower semicontinuous and $\Psi$ continuous and
satisfying $\lim_{x\in X, \|x\|\to +\infty}\Psi(x)=+\infty$. Put
$$\lambda^{*}=\inf_{\rho>\inf\limits _{X}\Psi}\inf_{x\in
\Psi^{-1}(]-\infty,\rho[)}
{{\Phi(x)-\inf_{\Psi^{-1}(]-\infty,\rho])}\Phi}\over
 {\rho-\Psi(x)}}\ .$$ 
Then, for each $\lambda>\lambda^{*}$, the functional $\Phi+\lambda
\Psi$ has a global minimum in $X$. Moreover, if $\lambda^{*}>0$,
for each $\mu<\lambda^{*}$, the functional $\Phi+\mu\Psi$
has no global minima in $X$.}\par
\bigskip
REMARK 1. - The second conclusion of Theorem 2 does not hold, in general,
if $\lambda^{*}=0$. To see this, consider, for instance, the case when
$\Phi(x)=\|x\|^2$ and $\Psi(x)=\|x\|$.\par
\bigskip
 Let us recall that a G\^ateaux differentiable functional $J$ on a real
 Banach
space $X$ is said to satisfy the Palais-Smale condition if
each sequence $\{x_{n}\}$ in $X$ such that $\sup_{n\in {\bf N}}|
J(x_{n})|<+\infty$ and $\lim_{n\to +\infty}\|J'(x_{n})\|_{X^{*}}=0$ admits
a strongly converging subsequence.\par
\bigskip
A joint application of conclusion $(a)$ of Theorem 1
with the main result of [13] gives the
following two critical points theorem:\par
\bigskip
THEOREM 3 ([17]). - {\it Let $X$ be a reflexive real Banach spaces and
let $\Phi, \Psi : X\to {\bf R}$ be two sequentially weakly lower
semicontinuous and continuously G\^ateaux differentiable functionals. 
 In addition, assume that, for each $\lambda>0$, the
functional $\Phi+\lambda\Psi$
 satisfies the Palais-Smale condition.\par
Then, for each $\rho>\inf_{X}\Psi$ and each $\mu$ satisfying 
$$\mu>\inf_{x\in \Psi^{-1}(]-\infty,\rho[)}
{{\Phi(x)-\inf_{\overline {(\Psi^{-1}(]-\infty,\rho[))}_{w}}\Phi}
\over {\rho-\Psi(x)}}\ ,$$
the following alternative holds: either the functional $\Phi+\mu
\Psi$ has a strict global minimum which lies
in $\Psi^{-1}(]-\infty,\rho[)$, or the same
functional has at least two critical points one of which lies
in $\Psi^{-1}(]-\infty,\rho[)$.}\par
\bigskip
We also recall that if $X$ is a real Hilbert space, an operator
$A:X\to X$ is said to be a potential operator if there exists a
a G\^ateaux differentiable functional $P$ on $X$ (which is called
a potential of $A$) such that $P'=A$.\par
\smallskip
Applying Theorem 1, we get the following result about fixed points of
potential operators in real Hilbert spaces:\par
\bigskip
THEOREM 4 ([14]). - {\it Let $X$ be a real Hilbert
space, and let $A:X\to X$ be a potential operator, with sequentially
weakly upper semicontinuous potential $P$. For each $\rho>0$,
put
$$\varphi(\rho)=\inf_{\|x\|^{2}<\rho}
{{\sup_{\|y\|^{2}\leq \rho}P(y)-P(x)}
\over {\rho-\|x\|^{2}}}\ .$$
Furthermore, set
$$\gamma=\liminf_{\rho\to +\infty}\varphi(\rho)$$
and
$$\delta=\liminf_{\rho\to {0}^{+}}\varphi(\rho)\ .$$
Then, the following conclusions hold:\par
\smallskip
\noindent
$(a)$\hskip 10pt If there is $\rho>0$ such that
$\varphi(\rho)<{{1}\over {2}}$,
then the operator $A$ has a fixed point whose norm is less than
$\sqrt {\rho}$.\par
\smallskip
\noindent
$(b)$\hskip 10pt If $\gamma<{{1}\over {2}}$, then
the following alternative holds: either the functional
$x\to {{1}\over {2}}\|x\|^{2}-P(x)$ has a global minimum, or the set
of all fixed points of $A$ is unbounded.\par
\smallskip
\noindent
$(c)$\hskip 10pt If $\delta<{{1}\over {2}}$, then
 the following alternative holds: either $0$ is 
a local minimum of the functional
$x\to {{1}\over {2}}\|x\|^{2}-P(x)$, or there exists a
sequence of pairwise distinct fixed points of $A$
which strongly converges to $0$.}\par
\bigskip
In particular, we have
\bigskip
THEOREM 5 ([14]). - {\it Let $X$ be a real Hilbert space, and let
$A:X\to X$ be a potential operator, with
 sequentially weakly upper semicontinuous 
potential $P$ satisfying 
$$\liminf_{r\to +\infty}{{\sup_{\|x\|\leq r}P(x)}\over {r^2}}<{{1}\over
{2}}<
\limsup_{r\to +\infty}{{\sup_{\|x\|\leq r}P(x)}\over {r^2}}\ .$$
Then, the set of all fixed points of $A$
is unbounded.}\par
\bigskip
>From now on, $\Omega$ is an open bounded subset of ${\bf R}^n$, with smooth
boundary, and (for $p>1$)
$W^{1,p}(\Omega)$, $W^{1,p}_{0}(\Omega)$ are the usual Sobolev
spaces, with norms
$$\|u\|=\left ( \int_{\Omega}|\nabla u(x)|^{p}dx + \int_{\Omega}|u(x)|^{p}dx
\right ) ^{1\over p}$$
and
$$\|u\|=\left ( \int_{\Omega}|u(x)|^{p}dx\right ) ^{1\over p}$$
respectively.\par
Let $p>1$, and  let $:\Omega\times {\bf R}\to {\bf R}$ be
 a Carath\'eodory function.\par
\smallskip
Recall that a weak solution of the Dirichlet problem
$$\cases {-\hbox {\rm div}(|\nabla u|^{p-2}\nabla u)
=f(x,u)
 & in
$\Omega$\cr & \cr u_{|\partial \Omega}=0\cr}$$
 is any $u\in W^{1,p}_{0}(\Omega)$ such that
 $$\int_{\Omega}|\nabla u(x)|^{p-2}\nabla u(x)\nabla v(x)dx
-\int_{\Omega}f(x,u(x))v(x)dx=0 \eqno{(*})$$
for all $v\in W^{1,p}_{0}(\Omega)$. While, a weak solution of the Neumann
problem
$$\cases {-\hbox {\rm div}(|\nabla u|^{p-2}\nabla u)=
f(x,u)
 & in
$\Omega$\cr & \cr {{\partial u}\over {\partial \nu}}=0 & on
$\partial \Omega$\cr} $$
 $\nu$ being the outer unit normal to $\partial \Omega$, is
 is any $u\in W^{1,p}(\Omega)$ satisfying identity $(*)$
for all $v\in W^{1,p}(\Omega)$. 
\bigskip
The following is one of the most classical existence results on the Dirichlet
problem for nonlinear elliptic equations:\par
\bigskip
THEOREM A ([1]). - {\it Assume that:\par
\noindent
$(1)$\hskip 5pt there are two positive constants $a, q$, with
$q<{{n+2}\over {n-2}}$ if $n\geq 3$, such that
$$|f(x,\xi)|\leq a(1+|\xi|^{q})$$
for all $(x,\xi)\in \Omega\times {\bf R}$;\par
\noindent
$(2)$\hskip 5pt there are constants $r\geq 0$ and $c>2$ such that
$$0<c\int_{0}^{\xi}f(x,t)dt\leq \xi f(x,\xi)$$
for all $(x,\xi)\in \Omega\times {\bf R}$ with $|\xi|\geq r$;\par
\noindent
$(3)$\hskip 5pt one has
$$\lim_{\xi\to 0}{{f(x,\xi)}\over {\xi}}=0$$   
uniformly with respect to $x$.\par
Then, the problem
$$\cases {-\Delta u=f(x,u)
 & in
$\Omega$\cr & \cr u_{|\partial \Omega}=0.\cr}$$
has a non zero weak solution.}\par
\medskip
What can be said whether, in Theorem A, condition $(3)$ is removed at all ?
Using Theorem 3, one obtains the following\par
\bigskip
THEOREM 6 ([17]). - {\it Assume that conditions $(1)$ and $(2)$ hold.\par
Then, for each $\rho>0$ and each $\mu$
satisfying
$$\mu>\inf_{u\in B_{\rho}}{{\sup_{v\in B_{\rho}}\int_{\Omega}
\left ( \int_{0}^{v(x)}f(x,\xi)d\xi\right )dx-
\int_{\Omega}
\left ( \int_{0}^{u(x)}f(x,\xi)d\xi\right )dx}\over
{\rho-\int_{\Omega}|\nabla u(x)|^{2}dx}}\ , \eqno{(**)}$$
where
$$B_{\rho}=\left \{ u\in W^{1,2}_{0}(\Omega):
\int_{\Omega}|\nabla u(x)|^{2}dx<\rho\right \}\  ,$$
the problem
$$\cases {-\Delta u={{1}\over {2\mu}} f(x,u)
 & in
$\Omega$\cr & \cr u_{|\partial \Omega}=0\cr}$$
has at least two weak solutions one of which lies in $B_{\rho}$.}\par
\bigskip
The following problem naturally arises in connection with Theorems A and
6.\par
\bigskip
PROBLEM 1. -  Under conditions $(1)$ and $(2)$ of Theorem A, is there
some $\rho>0$ such that the infimum appearing
in $(**)$ is less than ${{1}\over {2}}$ ? \par
\bigskip
Clearly, if the answer to this problem was positive, then Theorem 6 would be
a proper improvement of Theorem A.\par
\bigskip
Another result proved using Theorem 1 is the following:\par
\bigskip
THEOREM 7 ([18]). - {\it Let $f, g:\Omega\times {\bf R}\to {\bf R}$
be two Carath\'eodory functions. Assume that:\par
\noindent
$(i)$\hskip 10pt there is $s>1$ such that
$$\limsup_{\xi\to 0^+}
{{\sup_{x\in \Omega}|f(x,\xi)|}\over
{\xi^{s}}}<+\infty\ ;$$
\noindent
$(ii)$\hskip 10pt there is $q\in ]0,1[$ such that\par
$$\limsup_{\xi\to 0^+}{{\sup_{x\in \Omega}
|g(x,\xi)|}\over {\xi^q}}<+\infty\ ;$$
\noindent
$(iii)$\hskip 10pt there are a non-empty open set $D\subseteq \Omega$
and a set $B\subseteq D$ of positive measure such that
$$ \limsup_{\xi\to 0^{+}}{{\inf_{x\in B}
\int_{0}^{\xi}g(x,t)dt}\over {\xi^2}}=+\infty\ ,\hskip 3pt
 \liminf_{\xi\to 0^{+}}{{\inf_{x\in D}
\int_{0}^{\xi}g(x,t)dt}\over {\xi^2}}>-\infty\ .$$
Then, for some $\lambda^{*}>0$ and for
each $\lambda\in ]0,\lambda^{*}[$, 
the problem  
$$\cases {-\Delta u=
f(x,u)+\lambda g(x,u)
 & in
$\Omega$\cr & \cr u_{|\partial \Omega}=0\ , \cr} \eqno{\left ( P_{\lambda}
\right ) }$$
admits a non-zero, non-negative weak solution $u_{\lambda}\in
C^{1}(\overline {\Omega})$. Moreover, one has
$$\limsup_{\lambda\to 0^+}{{\|u_{\lambda}\|_{C^{1}(\overline {\Omega})}}\over
{\lambda^{q\over 1-q}}}<+\infty$$
and the function
$$\lambda\to {{1}\over {2}}\int_{\Omega}|\nabla u_{\lambda}(x)|^{2}dx-
\int_{\Omega}\left ( \int_{0}^{u_{\lambda}(x)}f(x,\xi)d\xi\right ) dx
-\lambda\int_{\Omega}\left ( \int_{0}^{u_{\lambda}(x)}g(x,\xi)d\xi
\right ) dx$$
 is negative
and decreasing in $]0,\lambda^{*}[$.
If, in addition,
 $f, g$ are continuous in $\Omega\times [0,+\infty[$
and
$$\liminf_{\xi\to 0^{+}}{{\inf_{x\in \Omega}g(x,\xi)}\over
{\xi |\log\xi|^{2}}}>-\infty\ ,$$
then $u_{\lambda}$ is positive in $\Omega$.}\par
\bigskip
REMARK 2. - Observe that Theorem 7 is bifurcation result. In fact, 
it ensures, in particular, that $\lambda=0$ is a bifurcation point
for problem $\left ( P_{\lambda}\right )$, in the sense that $(0,0)$
belongs to the closure in $C^{1}(\overline {\Omega})\times {\bf R}$
of the set
$$\left \{ (u,\lambda)\in C^{1}(\overline {\Omega})\times ]0,+\infty[
\hskip 3pt :\hskip 3pt
\hbox {\rm $u$ is a weak solution of $\left ( P_{\lambda}\right )$},
\hskip 3pt u\neq 0,\hskip 3pt u\geq 0\right \}\ .$$\par
\bigskip
In view of [2], where problem $\left ( P_{\lambda}\right )$ is studied for
particular nonlinearities, the following problems arises:\par
\bigskip
PROBLEM 2. - Under the assumptions of Theorem 7, does
problem $\left ( P_{\lambda}\right )$ admit a non-zero,
non-negative, {\it minimal} solution
for each $\lambda>0$ small enough ?    \par
\bigskip
REMARK 3. - Again applying Theorem 1, some bifurcations theorem for
Hammerstein nonlinear integral equations (in the spirit of Theorem 7) have
been obtained by F. Faraci in [10].\par
\bigskip
Note, in particular, the following corollary of Theorem 7:\par
\bigskip
COROLLARY 1 ([18]). - {\it Let $0<q<1<s$ and let
$\alpha, \beta$ be two H\"older continuous
functions on $\overline {\Omega}$. Assume that
$$0\leq \inf_{\Omega}\beta,\hskip 5pt 0<\sup_{\Omega}\beta\ . $$
Then, for some $\lambda^{*}>0$ and for each $\lambda\in ]0,\lambda^{*}[$,
the problem
$$\cases {-\Delta u=
\alpha(x) u^{s}+\lambda \beta(x) u^{q}
 & in
$\Omega$\cr & \cr
u_{|\partial \Omega}=0\cr}$$ 
admits a positive classical solution $u_{\lambda}$.
 Moreover, one has
$$\limsup_{\lambda\to 0^+}{{\|u_{\lambda}\|_{C^{1}(\overline {\Omega})}}\over
{\lambda^{q\over 1-q}}}<+\infty$$
and the function
$$\lambda\to {{1}\over {2}}\int_{\Omega}|\nabla u_{\lambda}(x)|^{2}dx-
{{1}\over {s+1}}\int_{\Omega}\alpha(x)|u_{\lambda}(x)|^{s+1}dx
-{{\lambda}\over {q+1}}\int_{\Omega}\beta(x)|u_{\lambda}(x)|^{q+1}dx$$
is negative and decreasing in $]0,\lambda^{*}[$.}\par
\bigskip
In the next theorems, $\lambda$ denotes a function in $L^{\infty}(\Omega)$,
with $\hbox {\rm ess inf}_{\Omega}\lambda>0$. They have been obtained in
[16] as applications of conclusions $(b)$ and $(c)$ of Theorem 1.\par
\bigskip
THEOREM 8 ([16]). - {\it Assume $p>n$.
Let $f:{\bf R}\to {\bf R}$ be a continuous function,
and $\{a_{k}\}$, $\{b_{k}\}$ two sequences
in ${\bf R}^+$ satisfying
$$a_{k}<b_{k}\hskip 5pt \forall k\in {\bf N},\hskip 5pt
\lim_{k\to \infty}b_{k}=+\infty,\hskip 5pt 
\lim_{k\to\infty}{{a_{k}}\over {b_{k}}}=0\ ,$$
$$\max\left \{ \sup_{\xi\in [a_{k},b_{k}]}\int_{a_{k}}^{\xi}
f(t)dt, \sup_{\xi\in [-b_{k},-a_{k}]}\int_{-a_{k}}^{\xi}
f(t)dt\right \}\leq 0\hskip 5pt \forall k\in {\bf N}$$
and
$$\limsup_{|\xi|\to +\infty}{{\int_{0}^{\xi}f(t)dt}\over
{|\xi|^{p}}}=+\infty\ .$$
Then, for every $\alpha, \beta\in L^{1}(\Omega)$, with $\min\{\alpha(x),
\beta(x)\}\geq 0$ a.e. in $\Omega$ and $\alpha\neq 0$,
and for every continuous function $g:{\bf R}\to {\bf R}$ satisfying
$$\sup_{\xi\in {\bf R}}\int_{0}^{\xi}g(t)dt\leq 0$$
and
$$\liminf_{|\xi|\to +\infty}{{\int_{0}^{\xi}g(t)}\over {|\xi|^{p}}}>
-\infty\ ,$$
the problem
$$\cases {-\hbox {\rm div}(|\nabla u|^{p-2}\nabla u)+\lambda(x) |u|^{p-2}u=
\alpha(x)f(u)+\beta(x)g(u) & in
$\Omega$\cr & \cr {{\partial u}\over {\partial \nu}}=0 & on
$\partial \Omega$\cr} $$
admits an unbounded sequence of weak solutions in $W^{1,p}(\Omega)$.}
\bigskip
A typical example of application of Theorem 8 is as follows:\par
\bigskip
EXAMPLE 1 ([16]). - Let $p>n$. Then, for each $\eta\in L^{1}(\Omega)$ with
 $\hbox {\rm ess inf}_{\Omega}\eta>0$, the problem
$$\cases {-\hbox {\rm div}(|\nabla u|^{p-2}\nabla u)
=\eta(x)\left ( \sum_{k=1}^{\infty}\left ( \hbox {\rm dist}
(u,{\bf R}\setminus [k!k,(k+1)!])\right ) ^{p}  
 -|u|^{p-2}u\right ) & in
$\Omega$\cr & \cr {{\partial u}\over {\partial \nu}}=0 & on
$\partial \Omega$\cr} $$
admits an unbounded sequence of weak solutions in $W^{1,p}(\Omega)$.\par
\bigskip
THEOREM 9 ([16]). - {\it Assume $p>n$.
Let $f:{\bf R}\to {\bf R}$ be a continuous function,
and $\{a_{k}\}$, $\{b_{k}\}$ two sequences
in ${\bf R}^+$ satisfying
$$a_{k}<b_{k}\hskip 5pt \forall k\in {\bf N},\hskip 5pt
\lim_{k\to \infty}b_{k}=0,\hskip 5pt 
\lim_{k\to\infty}{{a_{k}}\over {b_{k}}}=0\ ,$$
$$\max\left \{ \sup_{\xi\in [a_{k},b_{k}]}\int_{a_{k}}^{\xi}
f(t)dt, \sup_{\xi\in [-b_{k},-a_{k}]}\int_{-a_{k}}^{\xi}
f(t)dt\right \}\leq 0\hskip 5pt \forall k\in {\bf N}$$
and
$$\limsup_{\xi\to 0}{{\int_{0}^{\xi}f(t)dt}\over
{|\xi|^{p}}}=+\infty\ .$$
Then, for every $\alpha, \beta\in L^{1}(\Omega)$, with $\min\{\alpha(x),
\beta(x)\}\geq 0$ a.e. in $\Omega$ and $\alpha\neq 0$,
and for every continuous function $g:{\bf R}\to {\bf R}$ satisfying
$$\sup_{\xi\in {\bf R}}\int_{0}^{\xi}g(t)dt\leq 0$$
and
$$\liminf_{\xi\to 0}{{\int_{0}^{\xi}g(t)}\over {|\xi|^{p}}}>
-\infty\ ,$$
the problem
$$\cases {-\hbox {\rm div}(|\nabla u|^{p-2}\nabla u)+\lambda(x) |u|^{p-2}u=
\alpha(x)f(u)+\beta(x)g(u) & in
$\Omega$\cr & \cr {{\partial u}\over {\partial \nu}}=0 & on
$\partial \Omega$\cr} $$
admits a sequence of non zero weak solutions which strongly
converges to $0$ in $W^{1,p}(\Omega)$.}\par
\bigskip
REMARK 4. - Again using Theorem 1,
Theorems 8 and 9 have been extended to the setting of elliptic
variational-hemivariational inequalities by D. Motreanu and S. A. Marano in
[12].
While P. Candito, in [8], extended them to the case of
discontinuous nonlinearities. Further various applications of Theorem 1
can be found in [3], [4], [5], [9], [11] and [15]. See also the related papers
[6] and [19].
\par
\bigskip
 PROBLEM 3. - In Theorems 8 and 9, when $g=0$, are the conclusions still
valid without the assumption
$$\lim_{k\to\infty}{{a_{k}}\over {b_{k}}}=0\hskip 3pt ?$$ 
A partial answer to this problem has recently been provided by G. Anello and
G. Cordaro in [7].\par
\hfill\eject
\centerline {{\bf References}}\par
\bigskip
\bigskip
\noindent
[1]\hskip 5pt A. AMBROSETTI and P. H. RABINOWITZ, {\it Dual variational
methods in critical point theory and applications}, J. Funct. Anal.,
{\bf 14} (1973), 349-381.\par
\smallskip
\noindent
[2]\hskip 5pt A. AMBROSETTI, H. BREZIS and G. CERAMI, {\it Combined
effects of concave and convex nonlinearities in some elliptic problems},
J. Funct. Anal., {\bf 122} (1994), 519-543.\par
\smallskip
\noindent
[3]\hskip 5pt G. ANELLO and G. CORDARO, {\it Existence of solutions
of the Neumann problem involving the p-Laplacian via a variational
principle of Ricceri}, Arch. Math. (Basel), to appear.\par
\smallskip
\noindent
[4]\hskip 5pt G. ANELLO and G. CORDARO, {\it An existence theorem for the
Neumann problem involving the $p$-Laplacian}, preprint.\par
\smallskip
\noindent
[5]\hskip 5pt G. ANELLO and G. CORDARO, {\it An existence and localization
theorem of solution for a Dirichlet problem}, preprint.\par
\smallskip
\noindent
[6]\hskip 5pt G. ANELLO and G. CORDARO, {\it Positive infinitely many and
arbitrarily small solutions for the Dirichlet problem involving the
$p$-Laplacian}, Proc. Royal Soc. Edinburgh Sect. A, to appear.\par
\smallskip
\noindent
[7]\hskip 5pt G. ANELLO and G. CORDARO. {\it Infinitely many positive
solutions for the Neumann problem involving the $p$-Laplacian}, preprint.
\par
\smallskip
\noindent
[8]\hskip 5pt P. CANDITO, {\it Infinitely many solutions to the Neumann
problem for elliptic equations involving the $p$-Laplacian and with
discontinuous nonlinearities}, Proc. Math. Soc. Edinburgh, to appear.\par
\smallskip
\noindent
[9]\hskip 5pt G. CORDARO, {\it Existence and location of periodic solutions
to convex and non coercive Hamiltonian systems}, preprint.\par
\smallskip
\noindent
[10]\hskip 5pt F. FARACI, {\it Bifurcation theorems for
Hammerstein nonlinear integral equations}, Glasgow Math. J., to appear.\par
\smallskip
\noindent
[11]\hskip 5pt F. FARACI, {\it Multiplicity results for a Neumann problem
involving the $p$-Laplacian}, preprint.\par
\smallskip
\noindent
[12]\hskip 5pt S. A. MARANO and D. MOTREANU, {\it Infinitely many critical
points of non-differentiable functions and applications to a Neumann type
problem involving the $p$-Laplacian}, J. Differential Equations, to appear.
\par
\smallskip
\noindent
[13]\hskip 5pt P. PUCCI and J. SERRIN, {\it A mountain pass
theorem}, J. Differential Equations, {\bf 60} (1985), 142-149.\par
\smallskip
\noindent
[14]\hskip 5pt B. RICCERI, {\it A general variational principle and
some of its applications}, J. Comput. Appl. Math., {\bf 113}
(2000), 401-410.\par
\smallskip
\noindent
[15]\hskip 5pt B. RICCERI, {\it Existence and location of solutions to
the Dirichlet problem for a class of nonlinear elliptic equations},
 Appl. Math. Lett., {\bf 14} (2001), 143-148.\par
\smallskip
\noindent
[16]\hskip 5pt B. RICCERI, {\it Infinitely many solutions of the Neumann
problem for elliptic equations involving the $p$-Laplacian},
 Bull. London Math. Soc., {\bf 33} (2001), 331-340.\par
\smallskip
\noindent
[17]\hskip 5pt B. RICCERI, {\it On a classical existence
theorem for nonlinear elliptic equations}, in ``Experimental,
constructive and nonlinear analysis'', M. Th\'era ed., 275-278, CMS Conf.
Proc. {\bf 27}, Canad. Math. Soc., 2000.\par
\smallskip
\noindent
[18]\hskip 5pt B. RICCERI, {\it A bifurcation theorem for nonlinear
elliptic equations}, preprint.\par
\smallskip
\noindent
[19]\hskip 5pt J. SAINT RAYMOND, {\it On the multiplicity of the solutions
of the equation $-\Delta u=\lambda f(u)$}, J. Differential Equations,
to appear.\par
\bigskip
\bigskip
Department of Mathematics\par
University of Catania\par
Viale A. Doria 6\par
95125 Catania, Italy
\bye